\documentclass[11pt,a4paper]{article}
\usepackage{amsfonts}
\usepackage{eurosym}
\usepackage{amssymb}
\usepackage{amsmath}
\usepackage[linkcolor=blue,colorlinks=true,urlcolor=blue,citecolor=blue]{hyperref}
\usepackage{float}
\usepackage{array}

\newtheorem{theorem}{Theorem}

\newtheorem{example}[theorem]{Example}

\newtheorem{lemma}[theorem]{Lemma}

\newtheorem{problem}[theorem]{Problem}

\def\ds{\displaystyle}

\Large

\begin{document}

	\title{New Development of Homotopy Analysis Method for a Non-linear Integro-Differential Equations with initial conditions}
    \author{Z. K. Eshkuvatov$^{1,2}$\\
	$^{1}$Faculty of Ocean Engineering Technology and Informatics, University
	Malaysia \\
	Terengganu (UMT), Kuala Terengganu, Malaysia\\
	$^{2}$Independent reseacher, Faculty of Applied Mathematics and
	Intellectual Technology, \\
	National University of Uzbekistan (NUUz), Tashkent, Uzbekistan\\
	Corresponding author: $^{*}$E-mail: zainidin@umt.edu.my \\
    }
\date{}
\maketitle

\begin{abstract}
   Homotopy analysis method (HAM) was proposed by Liao in 1992 in his PhD thesis for non-linear problems and was applied it in many different problems of mathematical-physics and engineering. In this note, a new development of homotopy analysis method (ND-HAM) is demonstrated for non-linear integro-differential equation (NIDEs) with initial conditions. Practical investigations revealed that ND-HAM leads easy way how to find initial guess and it approaches to the exact solution faster than the standard HAM, modified HAM (MHAM), new modified of HAM (mHAM) and more general method of HAM (q-HAM). Two examples are illustrated to show the accuracy and validity of the proposed method. Five methods are compared in each example.
   \newline
     
     
   \noindent \textbf{Mathematics Subject Classification (AMS 2010)}: 34A08,
     26A33.
\end{abstract}

\section{Introduction}

   It is known that non-linear phenomena appear in many applications in scientific fields, such as fluid dynamics, solid-state physics, plasma physics, mathematical biology, chemical kinetics, aerodynamics, kinetic theory of gases, quantum mechanics, mathematical economics, and theory elasticity. Many problems in the fields of physics, engineering and science are modelled by non-linear integro-differential equations (IDEs) and the analytic solution of the non-linear problems is not known or can be found in rear cases. The obtained linear and non-linear integro-differential equations can be solved through many numerical methods. The existence of solution and approximation of the problems have been investigated by many authors viz: homotopy analysis method (HAM) \cite{Liao1, Liao2, Liao3, Liao4}, modified HAM \cite{MHAM1, MHAM2, MHAM4, MHAM5, MHAM6}, optimal HAM \cite{OHAM1, OHAM2, OHAM3}, q-HAM  \cite{q-HAM1, q-HAM2, Oq-HAM1, Oq-HAM3, Oq-HAM4},  for variety of problems. A number of series numerical methods were derived for integro-differential equations, such as variational iteration method \cite{Batiha1}, Taylor-successive approximation method \cite{NIDE1}, the Adomian decomposition method (ADM) \cite{NIDE2, NIDE3, NIDE4}, differential transform method \cite{NIDE5}, Wavelet-Galerkin method \cite{Avudainayagam}, the Tau method \cite{NIDE6}, modified HPM \cite{Zain1}, Laplace transform ADM \cite{Waleed}, finite element method \cite{Aloev} for partial differential equations (PDEs).
    
   \noindent In this note, the following IDEs are considered: 
   \begin{problem}
    	Non-linear Volterra-Fredholm integro-differential equations (VFIDEs) of order $p$
    	\begin{eqnarray}  \label{eqn_a1}
    	&& \ds u^{(p)}(t) + \sum_{j=1}^{p-1}a_j(t)u^{(j)}(t)=f(t)+\lambda_1
    	\int\limits_{a}^{t}K_1(t,s)F_1(u(s))ds \nonumber \\ 
    	&& \qquad\quad + \lambda_2
    	\int\limits_{a}^{b}K_2(t,s)F_2(u(s))ds,
    	\end{eqnarray}
    	with initial conditions 
    	\begin{equation} \label{eqn_a2}
    	u^{(k)}(a)=\alpha_k,\text{\ \ } k=0,...,p-1, \, p\in \mathbb{N}, \, (p\geq 2), 
    	\end{equation}
    	where $ t\in \Omega =\left[a, b\right]$ and $K_1,K_2:\Omega \times \Omega
    	\longrightarrow \mathbb{R}, \, f:\Omega \longrightarrow \mathbb{R}$ and $a_{j}: \Omega \longrightarrow \mathbb{R}, \, j=1,\cdots p-1 $ are known functions, $\lambda_1, \, \lambda_2$ are parameters and $F_1, F_2:C(\Omega ,\mathbb{R})\longrightarrow \mathbb{R}$ are non-linear function as well as $u(t)$ is unknown function to be determined
    \end{problem}
    For the non-linear IDEs \eqref{eqn_a1}-\eqref{eqn_a2}, application of ND-HAM (new development of HAM), standard HAM proposed by Liao \cite{Liao1}, modified HAM named (MHAM) developed by Bataineh et al. \cite{MHAM1}, new development of HAM (mHAM) initiated by Ayati et al. \cite{MHAM2}, and more general method of homotopy analysis method called (q-HAM) proposed by El-Tawil and Huseen \cite{q-HAM1} are described in details. 
    
    In Example 1 (Table 1 and 2), advantages of ND-HAM over standard HAM \cite{Liao1}, modified MHAM \cite{MHAM1}, mHAM \cite{MHAM2} and q-HAM \cite{q-HAM1}, are shown practically. Matlab codes are developed for five methods (ND-HAM, HAM, MHAM, mHAM, q-HAM) and shown that by suitable choice of $x_i(s)$ of $f(s)=x_0(s)+x_1(s)+\cdots +x_n(s)$, ND-HAM can give better results than the standard HAM and all modified HAM.
        
    Organization of the structure is as follows. In section 2, basic concept of HAM is recalled, few modified HAM (MHAM, mHAM, q-HAM) are described and the detailed description of the proposed method is presented. Section 3, describes detailed application of the standard HAM and ND-HAM for the Problem 1. In Section 4, three examples are solved to illustrate the performance of proposed method and comparisons with other methods are demonstrated. Finally, the paper ends with a conclusion and acknowledgement in Section 5.

\section{Basic idea of standard HAM, modified HAM (MHAM, mHAM, q-HAM) and ND-HAM}
    
    Basic idea of HAM is as follows. Let non-linear equation be given by 
    \begin{equation}  \label{eqb1}
    N[u(t)]=0.
    \end{equation}
    Liao \cite{Liao1} has constructed the zero-order deformation equation in the form
    \begin{equation}  \label{eqb2}
    (1-q)\pounds [\phi(t;q)-u_{0}(t)]=q\hbar H(t) \Big[N[\phi (t;q)]\Big], 
    \end{equation}
    where $\pounds $ is the linear operator,$\ q\in \lbrack 0,1$] is the
    embedding parameter, $\hbar \neq 0$ is an auxiliary parameter, $H(t)$ is auxilary function,  $u_{0}(t)$ is an initial guess of the solution $u(t)$ satisfying initial or boundary conditions and $\phi (t;q)$ is an unknown function to be determined depending on the variables $t$ and $q$ and satisfies the following equation
    \begin{equation}  \label{eqb3}
    \phi^{(i)}(t;0)=u_{0}^{(i)}(t), i=0,1,2,...
    \end{equation}
    When the parameter $q$ increases from $0$ to $1$, the homotopy solution $%
    \phi (t;q)$ varies from $u_{0}(t)$ to solution $u(t)$ of the original
    equation \eqref{eqb1}. Using the parameter $q$ as dummy variable, the
    function $\phi (t;q)$ can be expanded in Taylor series
    \begin{equation} \label{eqb4}
    \phi (t;q)=u_{0}(t)+\sum\limits_{m=1}^{+\infty }u_{m}(t)q^{m}, 
    \end{equation}
    where 
    \begin{equation*} \label{eqb5}
    u_{m}(t)=\frac{1}{m!}\left. \frac{\partial ^{m}\phi (t;q)}{\partial q^{m}}%
    \right\vert_{q=0}.
    \end{equation*}
    Assuming that the auxiliary parameter $\hbar $ in \eqref{eqb2} is properly selected so that the series \eqref{eqb4} is convergent when $q=1$ 
    \begin{equation} \label{eqb6}
    u(t)=u_{0}(t)+\sum\limits_{m=1}^{+\infty }u_{m}(t).
    \end{equation}
    Approximate solution of Eq. \eqref{eqb1} can be written as 
    \begin{equation} \label{eqb7}
    u(t)=u_{0}(t)+\sum\limits_{m=1}^{n}u_{m}(t),
    \end{equation}
    where $u_{m}(t)$ is defined from deformation equation \eqref{eqb9}.
    
    \noindent \textbf{High-order deformation equation}. We define the following
    vector $\overrightarrow{u}_{n}(t)=(u_{0}(t), \cdots ,u_{n}(t))$ and Differentiating the zero-order deformation equation \eqref{eqb2} $m$ times with respect to the embedding parameter $q$ and then dividing by $m!$ and finally setting $q=0,$ we have the so-called $m$th order deformation equation
    %
    \begin{equation} \label{eqb9}
    \pounds [u_{m}(t)-\chi _{m}u_{m-1}(t)]=\hbar H(t)
    \Re_{m}(\overrightarrow{u}_{m-1}(t)),  
    \end{equation}
    where 
    \begin{equation} \label{eqb10}
    \Re_{m}(\overrightarrow{u}_{m-1}(t))=\frac{1}{(m-1)!}\left. 
    \frac{\partial^{m-1}\Big[N[\phi (t;q)]\Big]}{\partial q^{m-1}} \right\vert_{q=0},  
    \end{equation}
    and 
    \begin{equation} \label{eqb11}
    \chi _{m}=\left\{ 
    \begin{array}{c}
    0,\text{ \ }m\leqslant 1, \\ 
    1,\text{ \ }m>1.%
    \end{array}%
    \right.  
    \end{equation}
    
    \noindent Series solution which are defined by \eqref{eqb9}-\eqref{eqb11} is called standard HAM. There are few modifications of HAM for instance: 
    \begin{itemize}
    	\item In 2009, Bataineh et al. \cite{MHAM1} had presented a new modification of the homotopy analysis method (MHAM) for solving non-linear systems of second order boundary-value problems (BVPs). In this modification, they consider the non-linear systems of second order differential equation in the form
        \begin{equation}  \label{eqb12}
        N_i[u_i(t)]=f_i(t), \, i=\{1,2\},
        \end{equation}
        subject to the boundary conditions:
        \begin{equation}  \label{eqb13}
        u_1(0)=u_1(1)=u_2(0)=u_2(1)=0.
        \end{equation}
        Representing the right side function 
        \begin{equation} \label{eqb14}
        f_i(t)= x_{i,0}(t)+ x_{i,1}(t)+ \cdots x_{i,n}(t), \, i=\{1,2\},
        \end{equation}
        and establishing $\varphi_i(t,q)$ in powers of the embedding parameter $q$ 
        \begin{equation} \label{eqb15}
        \varphi_i(t,q)= x_{i,0}(t)+ qx_{i,1}(t)+ \cdots +q^nx_{i,n}(t), \, i=\{1,2\},
        \end{equation}
        then constructing the $m$th-order deformation equation in the form
        \begin{eqnarray}  \label{eqb16}
        && \pounds [u_{i,m}(t)-\chi_{m}u_{i,m-1}(t)]=\hbar H(t)
        \Re_{i,m}(\overrightarrow{u}_{i,m-1}(t)), 
        \end{eqnarray}
        where $\chi_{m}$ is defined by \eqref{eqb11} and
        \begin{equation} \label{eqb17}
        \Re_{i,m}(\overrightarrow{u}_{i,m-1}(t))= \left. \frac{1}{(m-1)!} 
        \frac{\partial^{m-1}\Big[N_i[\phi_i(t;q)]-\varphi_i(t,q)\Big]}
        {\partial q^{m-1}}\right\vert_{q=0},
        \end{equation}
        is called a new modified version of HAM named MHAM.
        
        \item In 2014, Ayati et al. \cite{MHAM2} developed a new modified version of HAM named (mHAM). They have considered the following non-linear differential equation 
        \begin{equation}  \label{eqb18}
        N[u(t)]=f(t).
        \end{equation}
        The modified form of HAM, can be established based on assumption that the
        function $f(r)$ can be divided into several parts, namely,
        \begin{equation} \label{eqb19}
        f(t)= \sum_{k=0}^{n}x_k(t).
        \end{equation}
        Then, they have constructed the modified $m$th order deformation equation in the form
        \begin{eqnarray}  \label{eqb20}
        && \pounds [u_{0}(t)]=x_0(t),  \nonumber \\
        && \pounds [u_{1}(t) - u_0(t)]= \hbar\Big[\Re_{1}(\overrightarrow{u}_{0}(t)) -x_1(t) \Big],  \nonumber \\
        && \pounds [u_{m}(t) - u_{m-1}(t)]= \hbar\Big[\Re_{m}(\overrightarrow{u}_{m-1}(t)) - x_m(t) \Big], \, 2\leq m\leq n, \nonumber \\
        && \pounds [u_{m}(t)-u_{m-1}(t)]=\hbar\Re_{m}(\overrightarrow{u}_{m-1}(t)),  \, m>n.
        \end{eqnarray}
        \begin{equation} \label{eqb21}
        \Re_{m}(\overrightarrow{u}_{m-1}(t))=\frac{1}{(m-1)!}\left. 
        \frac{\partial^{m-1}\Big[N[\phi (t;q)]\Big]}{\partial q^{m-1}} \right\vert_{q=0},  
        \end{equation}
        where non-linear term $N[\phi (t;q)]$ is the left side of the equation \eqref{eqb18}.
         
        \noindent The series solution $u_m$ defined by \eqref{eqb19}-\eqref{eqb21} is called new modified HAM in short mHAM.  
        
        \item In 2012, El-Tawil and Huseen \cite{q-HAM1} considered the following differential equation:
        \begin{equation} \label{eqb211}
        N(u(x,t)) - f(x,t)= 0,
        \end{equation}
        where $N$ is a non-linear operator, $(x,t)$ denotes independent variables, $f(x,t)$ is a known function and $u(x,t)$ is an unknown function.
        
        \noindent They constructed the so-called zero-order deformation equation in the form
        \begin{equation}  \label{eqb212}
        (1-qn)\pounds [\phi(x,t,q)-u_{0}(x,t)]=q\hbar H(x,t) \Big[N[\phi(x,t,q)-f(x,t)]\Big], 
        \end{equation}
        where $n\geq 1, \, q\in \left[0, \frac{1}{n}\right]$, denotes the so-called embedded parameter, $L$ is an auxiliary linear operator with the property $L[f]=0$ when $f=0$, and $\hbar \not=0$ is an auxiliary parameter,  $H(x,t)$ denotes a non-zero auxiliary function.
        
        \noindent It is obvious, that when $q=0$ and $q=\frac{1}{n}$ equation \eqref{eqb212} becomes
        \begin{equation}  \label{eqb213}
        \phi(x,t,0) = u_{0}(x,t), \, \phi\left(x,t,\frac{1}{n}\right) = u(x,t), 
        \end{equation}
        respectively. Thus, as increases $q$ from $0$ to $\frac{1}{n}$, the solution $u(x,t)$ varies from the initial guess $u_0(x,t)$ to the solution $u(x,t)$.
        
        \noindent Expanding $\phi(x,t,q)$ in Taylor series with respect to $q$, one has
        \begin{equation} \label{eqb214}
        u(x,t)=u_{0}(x,t)+\sum\limits_{m=1}^{+\infty }u_{m}(x,t)q^m.
        \end{equation}
        where 
        \begin{equation*} \label{eqb215}
        u_{m}(x,t)=\frac{1}{m!}\left. \frac{\partial ^{m}\phi (x,t;q)}{\partial q^{m}}\right\vert_{q=0}.
        \end{equation*}
        Assume that $\hbar, \, H(x,t), \, u_0(x,t), \, L$ are so properly chosen such that the series \eqref{eqb214} converges at $q=\frac{1}{n}$ and the solution of \eqref{eqb211} can be written as
        \begin{equation} \label{eqb216}
        u(x,t)=\phi\left(x,t,\frac{1}{n}\right)= u_{0}(x,t) + \sum\limits_{m=1}^{+\infty }u_{m}(x,t)\left(\frac{1}{n}\right)^m.
        \end{equation}
        Differentiating equation \eqref{eqb212} $m$ times with respect to $q$ and  setting $q=0$ and finally dividing them by $(m-1)!$, we have the so-called $m-th$ order deformation equation:
        \begin{eqnarray}  \label{eqb217}
        && \pounds [u_{m}(x,t)-\chi_{m}u_{m-1}(x,t)]=\hbar H(x,t)
        \Re_{m}(\overrightarrow{u}_{m-1}(x,t)), 
        \end{eqnarray}
        where $\overrightarrow{u}_{r}(x,t)= (u_0(x,t), u_1(x,t), \cdots u_r(x,t))$ is a vector function and
        \begin{equation} \label{eqb218}
        \Re_{m}(\overrightarrow{u}_{m-1}(x,t))= \left. \frac{1}{(m-1)!} 
        \frac{\partial^{m-1}\Big[N[\phi(x,t,q)]-f(x,t)\Big]}
        {\partial q^{m-1}}\right\vert_{q=0},
        \end{equation}
        with
        \begin{equation} \label{eqb219}
        \chi _{m}=\left\{ 
        \begin{array}{c}
        0,\text{ \ }m\leqslant 1, \\ 
        n,\text{ \ }m>1.%
        \end{array}%
        \right.  
        \end{equation}
        \noindent The iterative solution $u_m(x,t)$ which are defined by \eqref{eqb217}- \eqref{eqb219} is called q-HAM. It is practically shown in El-Tawil and Huseen \cite{q-HAM1} that series solution \eqref{eqb216} where component functions $u_m(x,t)$ are governed by the equation \eqref{eqb217} converges to the exact solution $u(x,t)$, faster than the standard HAM. Due to the existence of factor $(1/n)^m$, more chances of convergence may occur. It should be noted that the case of $n=1$, standard HAM can be reached. In \cite{q-HAM2} convergence of q-HAM is proved and development of q-HAM and application of q-HAM are shown in \cite{Oq-HAM1, Oq-HAM3}.
        
        \item In 2015,  Yin et al. \cite{MHAM4} proposed a modified HAM to obtain quick and accurate solution of wave-like fractional physical models. This modified semi-analytical approach is the combination of Laplace transform algorithm and homotopy analysis method named homotopy analysis transform
        method (HATM) and was applied it for fractional partial differential equations. The HATM utilizes a simple and powerful method to adjust and control the convergence region of the infinite series solution using an auxiliary parameter. The numerical solutions obtained by this proposed method indicate that the approach is easy to implement, highly accurate, and computationally very attractive. A good agreement between the obtained solutions and some well-known results which have been obtained by other methods exist.
        
        \item In 2017, Ziane and Cherif \cite{MHAM5} proposed a new modification of HAM by combining the natural transformation with homotopy analysis method to solve non-linear fractional partial differential equations. This method is called the fractional homotopy analysis natural transform method (FHANTM). The FHANTM can easily be applied to many problems and is capable of reducing the size of computational work. The fractional derivative is described in the Caputo sense. The results show that the FHANTM is an appropriate method for solving non-linear fractional partial differential equation.
    \end{itemize}
    To derive the ND-HAM, we rewrite Eq. \eqref{eqb1} in the form 
    \begin{equation}  \label{eqb22}
    N[u(t)]=f(t),
    \end{equation}
    and assume that the function $f(t)$ is split into $n$ terms, namely
    \begin{equation} \label{eqb23} 
    f(t)=x_0(t)+x_1(t)+...+x_n(t).
    \end{equation}
    Expanding $g(t,q)$ into powers of the embedding parameter $q$, we obtain
    \begin{equation} \label{eqb24}
    g(t;q)=x_0(t)+x_1(t)+x_2(t)(qh)+...+x_n(t)(qh)^{n-1}.
    \end{equation}
    For ND-HAM we rewrite Eq. \eqref{eqb9}-\eqref{eqb11}, in the form
    \begin{eqnarray} 
    && \pounds [u_{0}(t)]=x_0(t), \label{eqb25} \\
    && \pounds [u_{m}(t)-\chi_{m}u_{m-1}(t)]=\hbar H(t)
    \Re_{m}(\overrightarrow{u}_{m-1}(t)), \label{eqb26}
    \end{eqnarray}
    where $\chi_{m}$ is defined by \eqref{eqb11} and
    \begin{equation} \label{eqb27}
    \Re_{m}(\overrightarrow{u}_{m-1}(t))= \left. \frac{1}{(m-1)!} 
    \frac{\partial^{m-1}\Big[N[\phi(t;q)]-g(t;q)\Big]}
    {\partial q^{m-1}}\right\vert_{q=0}.
    \end{equation}
    In ND-HAM, we have the following advantages: 
    \begin{enumerate}
    	\item Free choice of $x_0(t)$ depending on given function $f(t)$ and  solving Eq. \eqref{eqb25}, we can get exact solution $u(t)$ of the Eq. \eqref{eqb22} from the first iteration. In this case, the next iteration obtained by Eq. \eqref{eqb26} gives exactly zero solution $u_i(t) = 0, \, i=1,2,...$. 
    	
    	\item If Eq. \eqref{eqb25} does not give exact solution, then it will serve as a choice of initial guess $u_0(t)$ satisfying initial or boundary conditions. 
    	
    	\item Due to \eqref{eqb24}, residual term $\Re_{m}(\overrightarrow{u}_{m-1}(t))$ computed by \eqref{eqb27} gives us more economic computations in ND-HAM, for any value of $\hbar$ or equal computations compared to residual terms in HAM \eqref{eqb10}, MHAM \eqref{eqb17}, mHAM \eqref{eqb21} and q-HAM \eqref{eqb218}.
    	
    	\item Practical investigations reveal that ND-HAM might dominated standard HAM \cite{Liao1}, MHAM \cite{MHAM1}, mHAM \cite{MHAM2} standard HAM \cite{Liao1} and for small values of $n$ of Oq-HAM \cite{q-HAM1} (see Example 1).
    	
    \end{enumerate}

\section{\textbf{Main results. Application of HAM and ND-HAM}}
    
    To solve non-linear VFIDEs \eqref{eqn_a1}-\eqref{eqn_a2} using HAM and ND-HAM, we introduce a non-linear operator 
    \begin{eqnarray} \label{eqap1}
    && N\left[\phi (t;q)\right] =\Big(D_{a}^{p} +\sum_{j=1}^{p-1}a_{j}(t)D_{a}^{j}\Big)\phi(t;q) -\lambda_1 \int\limits_{a}^{t}K_1(t,s)F_1(\phi(s;q))ds \nonumber \\
    && \qquad\qquad\,\, - \lambda_2
    \int\limits_{a}^{b}K_2(t,s)F_2(\phi(s;q))ds,  
    \end{eqnarray}
    where $D^{j}_{a}$ is a differential operator of order $j$ and $\phi(t;q)$ is unknown function to be determined and write Eqs. \eqref{eqn_a1}-\eqref{eqn_a2} in the operator equation form
    \begin{eqnarray}
    && \ds N\left[\phi (t;q)\right] = f(t), \label{eqap2}\\
    && \phi^{(k)}(a;1) =u^{(k)}(a)= \alpha_k, \, k=0,1, \cdots p-1. \label{eqap3} 
    \end{eqnarray} 
    The $m^{\text{th}}$-order deformation equation (standard HAM) \eqref{eqb9}-\eqref{eqb11} with $H(t)=1$, for the non-linear equations \eqref{eqap2}-\eqref{eqap3} are given by
    \begin{eqnarray} \label{eqap4}
    && \displaystyle \pounds[u_{1}(t)]=\hbar \Re_{1}(\overrightarrow{u}_{0}(t)) = \hbar \left. \Big\{N[\phi(t;q)]-f(t)\Big\}\right\vert_{q=0}, \nonumber  \\
    && \displaystyle\pounds [u_{m}(t)- u_{m-1}(t)]= \hbar \Re_{m}(\overrightarrow{u}_{m-1}(t)),  \, m=2,3, \cdots \nonumber\\
    &&\displaystyle u_{0}^{(k)}(a)=\alpha_k, \,k=0,\cdots ,p-1, \, u_{m}^{(k)}(a)=0, \, m=1,2, \cdots, \, k=0, 1, \cdots  .
    \end{eqnarray}
    where residual term $\Re _{m}(\overrightarrow{u}_{m-1}(t))$ is as follows
    \begin{equation} \label{eqap5}
    \Re _{m}(\overrightarrow{u}_{m-1}(t))=\frac{1}{(m-1)!}\left[ \frac{\partial
    	^{m-1}[N(\phi (t;q))-f(t)]}{\partial q^{m-1}}\right]_{q=0}.  
    \end{equation}
    In a similar way, when right side function $f(t)$ of Eq. \eqref{eqap2} is split into the form \eqref{eqb23} - \eqref{eqb24}, then we can apply ND-HAM \eqref{eqb24}-\eqref{eqb27} for non-linear operator equation \eqref{eqap2}-\eqref{eqap3} as follows
    \begin{eqnarray} \label{eqap6}
    && \pounds [u_{0}(t)]=x_0(t), \nonumber\\
    && \pounds [u_{1}(t)]=\hbar \Re_{1}(\overrightarrow{u}_{0}(t)) = \left. \Big\{N[\phi(t;q)]-g(t,q)\Big\}\right\vert_{q=0},  \nonumber\\
    && \pounds [u_{m}(t)-u_{m-1}(t)]=\hbar \Re_{m}(\overrightarrow{u}_{m-1}(t)), \, m=2,3,\cdots \nonumber\\
    &&\displaystyle u_{0}^{(k)}(a)=\alpha_k, \, \, 
    k=0,\cdots, p-1, \, u_{m}^{(k)}(a)=0, \, m=1,2, \cdots, \, k=0,1, \cdots,
    \end{eqnarray}
    where $g(t,q)$ is defined by \eqref{eqb24} and residual term $\Re _{m}(\overrightarrow{u}_{m-1}(t))$ is of the form
    \begin{equation} \label{eqap7}
    \Re_{m}(\overrightarrow{u}_{m-1}(t))= \left. \frac{1}{(m-1)!} 
    \frac{\partial^{m-1}\Big[N[\phi(t;q)]-g(t,q)\Big]}
    {\partial q^{m-1}}\right\vert_{q=0}.
    \end{equation}
    Before applying HAM and ND-HAM for operator equation \eqref{eqap2}, we need the following important formulas Kanwal (\cite{Kanwal}, pp. 285)
    \begin{lemma}
    Let $f\in L^{1}[a,b]$ then $n$ tuple integrals $J^n_a(f)$ is computed in the form
    \begin{equation} \label{eqap8}
    J^n_a(f(s))=\int_{a}^{s}\int_{a}^{s_n}\cdots \int_{a}^{s_2}f(s_1)ds_1ds_2\cdots ds_n= \frac{1}{(n-1)!} \int_{a}^{s}(s-t)^{n-1}f(t)dt. 
    \end{equation}   	
    \end{lemma}
    \begin{lemma} 
    	Let $n, \, m \in N$ and $f\in L^{1}[a,b],$ then the the following properties hold
    \begin{eqnarray}
    	&&\ds J_{a}^{n}J_{a}^{m}f(t) =J_{a}^{m}J_{a}^{n}f(t) =J_{a}^{n+m} f(t), \notag \\
    	&&\ds D_{a}^{n}\left[J_{a}^{n}f(t) \right] = f(t),  \notag \\
    	&&\ds J_{a}^{n}\left[D_{a}^{n}f(t) \right] =f(t) -\sum_{k=0}^{n-1}\frac{f^{(k)}(a)}{k!}(t-a)^{k},  \label{eqap9} \\
    	&& \ds J_{a}^{m}(t-a)^{k}=\frac{k!}{(m+k)!}(t-a)^{m +k }.\,  \notag \\
    	&& \ds D_{a}^{n}\left(t-a\right)^{k}=0,\, k=0,1,...,n-1. \notag
    \end{eqnarray}
    \end{lemma}
    Let us apply ND-HAM \eqref{eqap6}-\eqref{eqap7} to the Eq. \eqref{eqn_a1} with initial conditions \eqref{eqn_a2}. Choosing $L=D_{a^{+}}^{p}$ and acting integral operator $J_{a^+}^{p}$ on both sides of first equation of \eqref{eqap6} taking into account \eqref{eqap3}, we have
    \begin{equation*} \label{eqap10}
    \ds J_{a}^{p}\left[D_{a}^{p}u_0(t)\right] = J_{a}^{p}(x_0(t)) = \frac{1}{(p-1)!} \int_{a}^{s}(s-t)^{p-1}x_0(t)dt,
    \end{equation*} 
    which leads to 
    \begin{eqnarray} \label{eqap11}
    	&& \ds u_0(t)= \sum_{k=0}^{p-1}\frac{\alpha_k}{k!}(t-a)^{k} + \frac{1}{(p-1)!} \int_{a}^{t}(t-s)^{p-1}x_0(t)ds.
    \end{eqnarray} 
    Since $u_{m}^{(k)}(a)=0, \, m=1,2, \cdots, \, k=0,1, \cdots$ and applying integral operator $J_{a}^{p}$ on both sides of the second equation  \eqref{eqap6} and taking $u_0(t)$ as an initial guess, we obtain
    \begin{eqnarray} \label{eqap12}
    && \ds u_1(t) = \hbar J_{a}^{p} \left[R_1(\overrightarrow{u_0}(t))\right] =\hbar J_{a}^{p}\Big\{N[\phi(t;q)]-g(t,q)\Big\}\vert_{q=0} \nonumber \\
    && \ds \qquad \, =\hbar J_{a}^{p} \left[\Big(D_{a}^{p} +\sum_{j=1}^{p-1}a_{j}(t)D_{a}^{j}\Big)u_0(t) - \Big(\lambda_1 G_1^0(t,q) + \lambda_2 G_2^0(t,q)\Big) - \Big(x_0(t)+x_1(t)\Big)\right], \nonumber \\
    && \ds \qquad \, =\hbar \Big[u_0(t) -\sum_{k=0}^{p-1}\frac{\alpha_k}{k!}(t-a)^{k}  + \sum_{j=1}^{p-1}J_{a}^{p}\Big(a_{j}(t)D_{a}^{j}u_0(t)\Big) \nonumber \\
    && \ds \qquad \, \ \ - J_{a}^{p}\Big(\lambda_1 G_1^0(t,q) + \lambda_2 G_2^0(t,q)\Big) - J_{a}^{p}\Big(x_0(t)+x_1(t)\Big)\Big],
    \end{eqnarray} 
    where
    \begin{equation} \label{eqap13}
    G_i^{0}(t,q)=\int\limits_{a}^{t}K_i(t,s)\Big[F_i(\phi(s;q))\Big]_{q=0}ds, \, i=1,2.
    \end{equation}
    For $2 \leq m \leq n$, due to \eqref{eqap1} and \eqref{eqap6}-\eqref{eqap7}, we obtain
    \begin{eqnarray} \label{eqap14}
    && \displaystyle\pounds [u_{m}(t)- u_{m-1}(t)]= \hbar \Re_{m}(\overrightarrow{u}_{m-1}(t)) =\frac{1}{(m-1)!}\left[\frac{\partial
    	^{m-1}[N[\phi (t;q)]-g(t,q)]}{\partial q^{m-1}}\right]_{q=0} \nonumber\\
    && \qquad \qquad = \Big[\Big(D^{p} +\sum_{j=1}^{p-1}a_{j}(t)D^{j}\Big) u_{m-1}(t) \nonumber \\ 
    && \qquad \qquad -\Big(\lambda_1 G_1^{m-1}(t,q) + \lambda_2 G_2^{m-1}(t,q)\Big) - x_{m}(t)\hbar^{m-1}\Big],
    \end{eqnarray}    
    where
    \begin{equation} \label{eqap15}
    G_i^{m-1}(t,q)\vert_{q=0}=\int\limits_{a}^{t}K_i(t,s)\frac{1}{(m-1)!}\Big[\frac{\partial^{m-1}}{\partial q^{m-1}}F_i(\phi(s;q))\Big]_{q=0}ds, \, i\in \{1,2\}.
    \end{equation}
    Similarly, for $m \geq n+1$ we get
    \begin{eqnarray} \label{eqap16}
    && \displaystyle\pounds [u_{m}(t)- u_{m-1}(t)]= \hbar \Re_{m}(\overrightarrow{u}_{m-1}(t)) =\frac{1}{(m-1)!}\left[\frac{\partial
    	^{m-1}[N[\phi (t;q)]]}{\partial q^{m-1}}\right]_{q=0} \nonumber\\
    && \qquad = \left[\Big(D^{p} +\sum_{j=1}^{p-1}a_{j}(t)D^{j}\Big) u_{m-1}(t) -\Big(\lambda_1 G_1^{m-1}(t,q)\vert_{q=0} + \lambda_2 G_2^{m-1}(t,q)\vert_{q=0}\Big) \right],
    \end{eqnarray}    
    where $G_i^{m-1}(t,q)\vert_{q=0}, \, i=1,2$ are defined by \eqref{eqap15}. 
     
    \noindent In view of $u_{m}^{(k)}(a)=0, \, m=1,2, \cdots, \, k=0,1, \cdots$ and acting integral operator $J_{a}^{p}$ on both sides of the equation  \eqref{eqap14} and taking into account the properties of integral operator $J_{a}^{p}$ \eqref{eqap8}-\eqref{eqap9}, we obtain
    \begin{eqnarray} \label{eqap17}
    && \displaystyle u_{m}(t) = (1+\hbar )u_{m-1}(t) + \sum_{j=1}^{p-1}J_{a}^{p}
    \Big(a_{j}(t)D_{a}^{j} u_{m-1}(t)\Big) \nonumber \\ 
    && \ds \qquad \quad - J_{a}^{p}\Big(\lambda_1 G_1^{m-1}(t,q) + \lambda_2 G_2^{m-1}(t,q)\Big) - J_{a}^{p}\Big(x_{m}(t)\Big)\hbar^{m-1}, \, 2\leq m \leq n,
    \end{eqnarray}    
    For $m \geq n+1$ and acting $J_{a}^{p}$ to both sides of Eq. \eqref{eqap16} together with the properties of integral operators \eqref{eqap8}-\eqref{eqap9}, we get
    \begin{eqnarray} \label{eqap18}
    && \displaystyle u_{m}(t) = (1+\hbar )u_{m-1}(t) + \sum_{j=1}^{p-1}J_{a}^{p}
    \Big(a_{j}(t)D_{a}^{j}\Big) u_{m-1}(t) \nonumber\\
    && \qquad \quad - J_{a}^{p}\Big(\lambda_1 G_1^{m-1}(t,q) + \lambda_2 G_2^{m-1}(t,q)\Big),
    \end{eqnarray}
    where $J_{a}^{p}f(t)$ is computed by \eqref{eqap8} and the derivative of
    the product of two functions $f(t)$ and $g(t)$ is computed by Leibniz rule, 
    \begin{equation}  \label{eqap19}
    D_{a}^{N}[f(t)g(t)]=\sum_{n=0}^{N}\binom{N}{n}D_{a}^{N-n}f(t)D_{a}^{n}g(t).
    \end{equation}

    \section{Numerical Experiments}
    
    In this section, we demonstrated five methods: Standard HAM \cite{Liao1}, MHAM \cite{MHAM1}, mHAM \cite{MHAM2} and q-HAM \cite{Oq-HAM1} in each of the example.
     
    \setcounter{theorem}{0}
    
    \begin{example} (Waleed Al-Hayani \cite{Waleed})
    	Let us solve the following second-order integro-differential equation  
   	\begin{equation} \label{eqex1} 
    	\left\{
    	\begin{array}{l} 
    	\ds y''(s)= e^{s} - s + \int_{0}^{1}sty(t)dt, \\ 
    	\ds y(0)=1,\, \, y'(0)=1.
    	\end{array}
    	\right.  
   	\end{equation} 
    	Exact solution is $y(s)=e^{s}$. 
    \end{example}

    \noindent {\bf Solution}: Let us find the exact and approximate  solution of problem \eqref{eqex1} using ND-HAM. Let us rewrite Eq. \eqref{eqex1} in the operator form
    \begin{equation}  \label{eqex2}
    \left\{ 
    \begin{array}{l}
    \ds N(\phi(s,q))= f(s), \\ 
    \phi(0,1) = y(0)=1, \, \phi'(0,1) = y'(0)=1, 
    \end{array}
    \right.
    \end{equation}
    where 
    \begin{equation} \label{eqex3}
    \left\{
    \begin{array}{l}
    \ds N(\phi(s,q))= \frac{\partial^2}{\partial s^2}\phi(s,q) - \int_{0}^{1}st\phi(t,q)dt, \\
    \ds f(s) = e^{s} - s.
    \end{array}
    \right.
    \end{equation}
    We attempt to find exact solution of IDEs \eqref{eqex1} by four methods.
    \begin{enumerate}
    \item {\bf ND-HAM}. Let us expand right side function of Eq. \eqref{eqex2} in the form  
    \begin{equation} \label{eqex4}
    f(s)= e^{s} - s = x_0(s) + x_1(s) = [g(s;q)]_{q=0},
    \end{equation}
    where $x_0(t) = e^{s}, \, x_1(s)= -s$. Choosing $\ds L=\frac{d^2}{ds^2}$, and from the first equation of \eqref{eqap6}, we have
    \begin{equation} \label{eqex5}
    L[y_0(s)]=x_0(s) \rightarrowtail \, y_0(s)= e^s.
    \end{equation}
    Let $m=1$, then from second equation of \eqref{eqap6}, we obtain
    \begin{eqnarray} \label{eqex6}
    && \ds \pounds [y_1(s)]=\hbar R_1(y_0(s))=h[N(\phi(s;q))-g(s;q)]\mid_{q=0} \nonumber\\
    && \ds \qquad \qquad =\hbar\Big[\frac{d^2}{ds^2}[y_0(s)]-s\int_0^{1}ty_0(t)dt -(x_0(s)+x_1(s))\Big] \nonumber \\
    && \qquad \qquad =\hbar\left[e^s -s - (e^s -s)\right] =0,
    \end{eqnarray}
    which implies $y_1(s)=0$. By continuing this procedure, we obtain $y_2(s)=y_3(s)=y_4(s)=\cdots =0$. So that the solution of equation \eqref{eqex2} is
    \begin{equation} \label{eqex7}
    y(s)=\phi(s,1)=\sum_{m=0}^{\infty} y_m(s) = y_0(t) = e^s,
    \end{equation}
    which is identical with the exact solution of Eq. \eqref{eqex1}. 

    \item {\bf Standard HAM}. Assume that initial guess is $y_0(s)=e^s$, then from   Eqs. \eqref{eqb9}-\eqref{eqb10} it follows that
    \begin{eqnarray} \label{eqex8}
    && \ds \pounds [y_1(s)]=\hbar R_1(y_0(s))=h[N(\phi(s;q))-f(s)]\mid_{q=0} \nonumber\\
    && \ds \qquad \qquad =\hbar\left[\frac{d^2}{ds^2}[y_0(s)]-s\int_0^{1}ty_0(t)dt -f(s)\right] \nonumber \\
    && \qquad \qquad =\hbar\left[e^s -s - (e^s -s)\right] =0.
    \end{eqnarray}
    Since $y_1(0)=y'_1(0)=0$ by integrating twice \eqref{eqex8}, we get $y_1(s)=0$. Continuing in this manner, we obtain $y_2(s)=y_3(s)=y_4(s)=\cdots =0$. Thus, the solution of the equation \eqref{eqex2} is
    \begin{equation} \label{eqex9}
    y(s)=\phi(s,1)=\sum_{m=0}^{\infty} y_m(s) = y_0(t) = e^s
    \end{equation}
    which is identical with the exact solution of Eq. \eqref{eqex1}. 
    
    \noindent Since first iteration of q-HAM developed by El-Tawil and Huseen \cite{q-HAM1} coincides with first iteration of HAM, therefore it gives exact solution when initial guess is chosen as exact solution.

   	\item {\bf mHAM} is developed by Ayati et al. \cite{MHAM2} in 2014.
   	In this case function $f(s)$ is divided into two parts 
    \begin{equation} \label{eqex10}
    f(s)=e^s-s= x_0(s)+x_1(s).
    \end{equation}
    where $x_0(s)=e^s, \, x_1(s)=-s$. 
    
    \noindent In view of the first equation of \eqref{eqb20}, we get
    \begin{equation} \label{eqex11}
    L[y_0(s)]=x_0(s) \rightarrowtail \, y_0(s)= e^s.
    \end{equation}
    From the second equation of \eqref{eqb20}, it follows that
    \begin{eqnarray} \label{eqex12}
    && \ds \pounds [y_1(s)-y_0(s)]=\hbar [R_1(y_0(s))-x_1(s)] \nonumber\\
    && \ds \qquad \qquad =\hbar\left[\frac{d^2}{ds^2}[y_0(s)]-s\int_0^{1}ty_0(t)dt -x_1(s)\right] =\hbar e^s,
    \end{eqnarray}
    which leads to 
    \begin{eqnarray} \label{eqex13}
    \ds y_1(s) =(1+\hbar)(e^s-1-s).
    \end{eqnarray}
    From the third equation of \eqref{eqb20}, we obtain second iteration
    \begin{eqnarray} \label{eqex14}
    && \ds \pounds [y_2(s)-y_1(s)]=\hbar R_2(y_1(s)) =\hbar\frac{\partial}{\partial q}\Big[N(\phi(s;q))\Big]_{q=0} \nonumber\\
    && \ds \qquad \qquad =\hbar\left[\frac{d^2}{ds^2}[y_1(s)]- s\int_0^{1}ty_1(t)dt\right] =\hbar(1+\hbar)\left[ e^s - \frac{s}{6} \right],
    \end{eqnarray}
    Acting operator $J^2_{0}$ on both sides of Eq. \eqref{eqex14} leads to
    \begin{eqnarray} \label{eqex15}
    && \ds y_2(s)= y_1(s) + \hbar(1+\hbar)\int_{0}^{s}\int_{0}^{s_1}\left[e^t-\frac{t}{6}\right]dtds_1 \nonumber \\
    && \ds \qquad\quad = (1+\hbar)(e^s-1-s) + \hbar(1+\hbar)\left(e^s-1-s-\frac{s^3}{36} \right) \nonumber \\
    && \ds \qquad\quad = (1+\hbar)^2(e^s-1-s) - \hbar (1+\hbar)\frac{s^3}{36}.
    \end{eqnarray}
    Three term iterations at $\hbar=-1$ yields
    \begin{equation} \label{eqex16}
    Y_{2mHAM}(s)= y_0(s)+ y_1(s)+ y_2(s) = e^s,
    \end{equation}
    which is identical with the exact solution of Eq. \eqref{eqex1}.

    \item {\bf In the modified HAM (MHAM)} developed by Bataineh \cite{MHAM1}, function $\varphi(t)$ is constructed as
    \begin{equation} \label{eqex17}
    \varphi(s,q)=x_0(s)+qx_1(s).
    \end{equation}
    where $\ds x_0(s)=e^s-\frac{4s}{5}, \, x_1(s)=-\frac{s}{5}$. 
    
    \noindent Choosing $y_0(s)=exp(s)$ and from \eqref{eqb16}-\eqref{eqb17} and \eqref{eqex17} it follows that
    \begin{eqnarray} \label{eqex18}
    && \ds \pounds [y_1(s)]=\hbar R_1(u_0(s))=h[N(\phi(s;q))-\varphi(s;q)]\vert_{q=0} \nonumber\\
    && \ds \qquad \qquad =\hbar\left[\frac{d^2}{ds^2}[y_0(s)]-s\int_0^{1}ty_0(t)dt -x_0(s)\right] =-\frac{4}{5}\hbar s,
    \end{eqnarray}
    which leads to 
    \begin{eqnarray} \label{eqex19}
    \ds y_1(s)=-\frac{4}{5}\hbar \int_{0}^{s}\int_{0}^{s_1}tdtds_1= -\hbar \frac{s^3}{30}.
    \end{eqnarray}
    Next iteration is
    \begin{eqnarray} \label{eqex20}
    && \ds \pounds [y_2(s)-y_1(s)]=\hbar R_2(y_1(s)) =\hbar\frac{\partial}{\partial q}\Big[N(\phi(s;q))-\varphi(s;q)\Big]_{q=0} \nonumber\\
    && \ds \qquad \qquad =\hbar\left[\frac{d^2}{ds^2}[y_1(s)]-s\int_0^{1}ty_1(t)dt -x_1(s)\right] =\hbar s\left[\frac{1}{5}-\hbar\frac{29}{150}\right],
    \end{eqnarray}
    Acting operator $J^2_{0}$ on both sides of Eq. \eqref{eqex20} leads to
    \begin{eqnarray} \label{eqex21}
    && \ds y_2(s)= y_1(s) + \hbar\left[\frac{1}{5}-\hbar\frac{29}{150}\right]
    \int_{0}^{s}\int_{0}^{s_1} tdtds_1 = -\hbar^2 \frac{29}{900}s^3.
    \end{eqnarray}
    Three term iterations at $\hbar=-1$ yields
    \begin{equation} \label{eqex22}
    Y_{2MHAM}(s)= y_0(s)+ y_1(s)+ y_2(s) = e^s + \frac{1}{900}s^3,
    \end{equation}
    which gives highly accurate solution but does not coincide with the exact solution.
    
    \end{enumerate}
    Let us apply ND-HAM for different initial guess. To find approximate solution let us split right side function of Eq. \eqref{eqex2} in the form
    \begin{equation} \label{eqex23}
    f(s)= 0 + \left(e^s - \frac{3s}{4}\right) - \frac{s}{4} = x_0(s) + x_1(s) + x_2(s),
    \end{equation}
    where $x_0(s) =0, \, x_1(s)=e^s-3s/4, \,  x_2(s)= -s/4$ and construct $g(s,t)$ function as follows
    \begin{equation} \label{eqex24}
    g(s,q) = x_0(s) + x_1(s) + (q\hbar)x_2(s).
    \end{equation}
    Knowing $\ds L=\frac{d^2}{ds^2}$, and solving the first equation of \eqref{eqap6}, we obtain
    \begin{equation} \label{eqex25}
    y_0(s) = 1 + s.
    \end{equation}
     In Table 1, we summarise all five methods with the same initial guess and with five iterations \\
     \begin{table}
     	\begin{flushleft} \centering \small
        \begin{tabular}{|c|c|c|c|c|}
    	\hline
    	$m$ & Exact & ND-HAM  & HAM & MHAM \\   \hline
    	3  & $\ds e^s$ & $\ds e^s+\frac{s^3}{2160}$ & $\ds e^s-\frac{s^3}{1080}$  & $\ds e^s+\frac{s^3}{540}$ \\ \hline
    	5  & $\ds e^s$ & $\ds e^s+\frac{s^3}{5832}\cdot 10^{-4}$ & $\ds e^s-\frac{s^3}{2916}\cdot 10^{-4}$ &  $\ds e^s+\frac{s^3}{1458}\cdot 10^{-4}$ \\ \hline
    	10 & $\ds e^s$ & $\ds e^s+\frac{s^3}{1417176}\cdot 10^{-9}$ & $\ds e^s-\frac{s^3}{708588}\cdot 10^{-9}$ & $\ds e^s+\frac{s^3}{354294}\cdot 10^{-9}$\\ \hline
        \end{tabular} 
   	    \caption{Numerical results for Example 1}
        \end{flushleft}
   \end{table}
   \begin{table}
       \begin{flushleft} \centering \small
       \begin{tabular}{|c|c|c|c|}
       \hline
    	$m$ & Exact & mHAM & q-HAM \\  \hline
   		3  & $\ds e^s$ & $\ds e^s+\frac{s^3}{540}$ & $\ds e^s-\frac{s^3}{1080}$ \\ \hline
   		5  & $\ds e^s$ & $\ds e^s+\frac{s^3}{1458}\cdot 10^{-4}$ & $\ds e^s-\frac{s^3}{2916}\cdot 10^{-4}$ \\ \hline
   		10 & $\ds e^s$ & $\ds e^s+\frac{s^3}{354294}\cdot 10^{-9}$ & $\ds e^s-\frac{s^3}{708588}\cdot 10^{-9}$\\ \hline
        \end{tabular} 
        \caption{Numerical results for Example 1}
        \end{flushleft}
     \end{table} 

    \noindent Here 
    \begin{itemize}
    	\item For MHAM function $\varphi(s,q)=x_0(s)+qx_1(s)$, where $\ds x_0(s)=e^s-\frac{s}{2},  \, x_1(s)=-\frac{s}{2} $
    	\item For mHAM function $f(s)=x_0(s)+x_1(s)+x_2(s)$, where $\ds x_0(s)=0, x_1(s)=e^s-\frac{s}{2},  \, x_2(s)=-\frac{s}{2} $
    	\item For q-HAM function $f(s)=e^s-s$ and $n=2, \hbar=-2$.
    \end{itemize}
 
    \noindent From the Table 1 and 2, we can conclude that ND-HAM dominated over all modified HAM. Control parameter $\hbar=-1$ is taken for ND-HAM, HAM, MHAM, mHAM and $\hbar=-2$ for q-HAM.
    
    \begin{example} (Huseen et al. \cite{Oq-HAM4})  Let us consider non-linear VIEs
   	\begin{eqnarray} \label{eqex281}
    	&& \ds u'(s) = -1+ \int_{0}^{s}u^2(t)dt, \nonumber \\
    	&& \ds u(0) = 0.
   	\end{eqnarray}
    \end{example}
    {\bf Solution}: Analytical solution of Eq. \eqref{eqex281} is not exist, fortunately Wavelet Galerkin method \cite{Avudainayagam} is taken as analytical numerical solution. To solve \eqref{eqex281} by the ND-HAM we rewrite it in the operator form
    \begin{equation}  \label{eqex282}
    \ds N(\phi(s,q))= f(s),
    \end{equation}
    where 
    \begin{equation} \label{eqex283}
    \left\{
    \begin{array}{l}
    \ds N(\phi(s,q))= \frac{\partial}{\partial s}\phi(s,q) -  \int_{0}^{s}\phi^2(t,q)dt, \\
    \ds f(s) = -1.
    \end{array}
    \right.
    \end{equation}
    Let us construct $g(s,q)$ function as follows
    \begin{equation} \label{eqex2831}
    \ds g(s,q)\vert_{q=0}=f(s)=-1 + 0 = x_0(s)+x_1(s), \, x_0(s)=-1, \, x_1(s)=0.
    \end{equation}
    In view of Eq. \eqref{eqb25} and \eqref{eqex2831}, it follows that $u_0(s)=-s$. Since 
    \begin{equation} \label{eqex284}
    \left\{
    \begin{array}{l} 
    \ds \phi(s,q) = u_0(s) + \sum_{m=1}^{\infty}q^m u_m(s), \\
    \ds \frac{\partial}{\partial q}\phi(s,q)\vert_{q=0} = u_1(s), \\
    \ds \frac{\partial^{m-1}}{\partial q^{m-1}}\phi(s,q)\vert_{q=0} = 
    (m-1)!u_{m-1}(s), \\
    \ds \phi^2(s,q)\vert_{q=0} = u_0^2(s), \, \frac{\partial}{\partial q}\phi^2(s,q)\vert_{q=0} = 2u_0(s)u_1(s),\\
    \ds \frac{1}{(m-1)!}\frac{\partial^{m-1}}{\partial q^{m-1}}\phi^2(s,q)\vert_{q=0} =  \sum_{i=0}^{m-1}u_i(s)u_{m-1-i}, \\
    \end{array}
    \right.    
    \end{equation}
    we have
    \begin{eqnarray} \label{eqex285}
    && \ds L[u_1(s)] = \hbar\left[N(\phi(s,q))-(x_0(s)+x_1(s))\right]\vert_{q=0} = \hbar\left[\frac{\partial}{\partial s}u_0(s) - \int_{0}^{s}u_0^2(t)dt -(x_0(s)+x_1(s))\right], \nonumber \\
    && \ds L[u_2(s) - u_1(s)] = \hbar\frac{\partial}{\partial q} \left[N(\phi(s,q))\right]\vert_{q=0} = \hbar \left[\frac{\partial}{\partial s}u_1(s) - 2\int_{0}^{s}u_0(t)u_1(t)dt \right], \nonumber\\
    && \ds L[u_m(s) - u_{m-1}(s)] = \hbar\frac{1}{(m-1)!}\frac{\partial^{m-1}}{\partial q^{m-1}} \left[N(\phi(s,q))\right]\vert_{q=0} \nonumber \\
    && \qquad \qquad = \hbar \left[\frac{\partial}{\partial s}u_{m-1}(s) - \int_{0}^{s}\sum_{i=0}^{m-1}u_i(t)u_{m-1-i}(t)dt \right], \, m\geq 3.
    \end{eqnarray}
    From \eqref{eqex285}, it follows that
    \begin{eqnarray} \label{eqex286}
    && \ds u_1(s) = -\frac{\hbar}{12} \cdot s^4,\nonumber \\
    && \ds u_2(s) = u_1(s) - \frac{\hbar^2}{252}s^4\cdot \left(21+s^3\right), \nonumber\\
    && \ds u_3(s) = u_{2}(s) + \frac{\hbar^2}{6048}s^4\left(504\cdot (1+h) + 24(1+2h)s^3 +hs^6\right).
    \end{eqnarray}
    Fifth terms approximation of the ND-HAM at $\hbar=-1$ is
    \begin{eqnarray} \label{eqex287}
    && U5NDHAM(s) = u_0(s) + u_1(s) + u_2(s) + u_3(s) + u_4(s) + u_5(s) \nonumber \\ 
    && \qquad \quad = -s+\frac{1}{12}s^4-\frac{1}{252}s^7 + \frac{1}{6048}s^{10}  - \frac{1}{157248}s^{13}+\frac{37}{158505984}s^{16}.
    \end{eqnarray}
     Since $f(s)=-1$ is constant, fifth terms approximation of the HAM, MHAM and mHAM are the same as ND-HAM defined by \eqref{eqex287}. Fifth terms approximation of the Oq-HAM, developed by Sh.N. Huseen et al. \cite{Oq-HAM4} in 2013, for the same problem of Eq. \eqref{eqex281}, at $n=2$ yields
    \begin{eqnarray} \label{eqex290}
    && U5OqHAM(s) = u_0(s) + \frac{u_1(s)}{2} + \frac{u_2(s)}{4} + \frac{u_3(s)}{8} + \frac{u_4(s)}{16} + \frac{u_5(s)}{32} \nonumber \\ 
    && \qquad = -s+\frac{1}{24}s^4-\frac{1}{1008}s^7 + \frac{1}{48384}s^{10}  - \frac{1}{2515968}s^{13}+\frac{37}{5072191488}s^{16}.
    \end{eqnarray}
    Fifth terms approximation of the Adomian decomposition method (ADM) developed in El-Sayed and Abdel-Aziz \cite{NIDE4} has the form,
    \begin{eqnarray} \label{eqex288}
    && U5ADM(s) = u_0(s) + u_1(s) + u_2(s) + u_3(s) + u_4(s) + u_5(s) \nonumber \\ 
    && \qquad \quad = -s+\frac{1}{12}s^4-\frac{1}{252}s^7 + \frac{1}{6048}s^{10} - \frac{1}{157248}s^{13}+\frac{79}{264176640}s^{16}.
    \end{eqnarray}

   \noindent Numerical comparisons of the methods are given in the Table 3. Five iterations are taken for all methods. Parameter $n=2$ is taken for the Oq-HAM \cite{Oq-HAM4}. \\

\begin{table}
    \begin{flushleft} \centering \small
    \begin{tabular}{|c|c|c|c|c|}
    	\hline
    	$s$ value & Exact \cite{Avudainayagam} & ND-HAM & Oq-HAM in \cite{Oq-HAM4}& ADM in \cite{NIDE4} \\
    	\hline
    	0.0000& 0.0000 & 0.000000000 & 0.0000000000 & 0.0000000 \\
    	\hline
    	0.0938&-0.0937 &-0.093793549 & -0.09379677  & -0.093793549 \\
    	\hline
    	0.3125&-0.3117 &-0.311706425 & -0.31210292  & -0.311706425 \\
    	\hline
    	0.5000&-0.4948 &-0.494822508 & -0.49740330  & -0.494822508  \\
    	\hline
    	0.7188&-0.6969 &-0.696941464 & -0.70776777  & -0.696941463 \\
    	\hline
    	0.9062&-0.8520 &-0.851934173 & -0.86852216  & -0.851934160 \\
    	\hline
    	1.000 &-0.9205 &-0.920475703 & -0.92911909  & -0.920475637 \\
    	\hline
   \end{tabular} \\
   \caption{Numerical results for Example 2}
   \end{flushleft}
\end{table}
    
    \noindent For the numerical comparisons of the methods Wavelet Galerkin method \cite{Avudainayagam} is taken as analytical numerical solution. From Table 3 we can note that ND-HAM which is identical with HAM is a bit more accurate than others. \\

\section{Conclusion}

In this paper, we explored the superiority of the ND-HAM over the HAM \cite{Liao1}, MHAM \cite{MHAM1}, mHAM \cite{MHAM2} and \cite{q-HAM1} in terms of accuracy. It provides a good opportunity to choose initial guess by solving linear operator equation depending on the given function $f(t)$. We also investigated the non-linear IDEs \eqref{eqap1}-\eqref{eqap2} and its scheme of approximate solutions. Five methods were compared with the same examples and the same initial conditions. Table 1-2 demonstrated that ND-HAM can give better results over other methods because of suitable choice $x_i(s)$ of the expansion $f(s)$. Example 2, are also solved by Oq-HAM \cite{Oq-HAM4} and Adomian decomposition method (ADM) \cite{NIDE3} until fifth iterations. Practical investigations showed that ND-HAM gave exact solution when initial guess coincides with the exact solutions. Even if it does not coincide with the exact solution, it brings one closer to the exact solutions with small number of iterations. Standard HAM is well developed and converges to exact solution fast and accurate. Matlab code is used to get the numerical solutions.\\

\noindent \textbf{Acknowledgement}. This work was supported by Universiti
Malaysia Terengganu (UMT) under RMC Research Grant Scheme (UMT, 2020).
Project code is UMT/CRIM/2-2/2/14 Jld. 4(44).

\end{document}